\documentclass[11pt]{article}
\usepackage{amssymb,latexsym,geometry,tikz,setspace}
\newtheorem{theorem}{Theorem}
\newtheorem{proposition}{Proposition}

\newtheorem{claim}{Claim}

\begin{document}

\onehalfspace

\title{Induced Matchings in Subcubic Graphs}

\author{Felix Joos, Dieter Rautenbach, Thomas Sasse}

\date{}

\maketitle

\vspace{-1cm}

\begin{center}
Institut f\"{u}r Optimierung und Operations Research, 
Universit\"{a}t Ulm, Ulm, Germany\\
\{\texttt{felix.joos, dieter.rautenbach, thomas.sasse}\}\texttt{@uni-ulm.de}
\end{center}

\begin{abstract}
We prove that a cubic graph with $m$ edges has an induced matching with at least $m/9$ edges.
Our result generalizes a result for planar graphs due to Kang, Mnich, and M\"{u}ller
(Induced matchings in subcubic planar graphs, SIAM J. Discrete Math. 26 (2012) 1383-1411)
and solves a conjecture of Henning and Rautenbach
(Induced matchings in subcubic graphs without short cycles, to appear in Discrete Math.).
\end{abstract}


\pagebreak

\section{Introduction}

We consider finite, simple, and undirected graphs and use standard terminology and notation.
An {\it induced matching} of a graph $G$ is a set $M$ of edges of $G$
such that every two vertices of $G$ 
that are incident with distinct edges in $M$
have distance at least $2$ in $G$.
The maximum number of edges in an induced matching of $G$ is the {\it strong matching number $\nu_s(G)$ of $G$}. 
The minimum number of induced matchings of $G$ 
into which the edge set of $G$ can be partitoned is the {\it strong chromatic index $\chi'_s(G)$ of $G$}.

The strong matching number and the strong chromatic index 
are analogs of the classical notions of the matching number $\nu(G)$ and the chromatic index $\chi'(G)$ of a graph $G$, respectively.
While matchings in graphs are very well understood \cite{lopl,tu} and the matching number can be determined efficiently \cite{ed},
Stockmeyer and Vazirani \cite{stva} proved the computational hardness of the strong matching number.
Their result was strengthened in many ways and restricted graph classes where the strong matching number
can be determined efficiently were studied \cite{brmo,ca1,ca2,gohehela,lo}.
Similarly, the chromatic index of a graph is much better understood than its strong variant.
While $\chi'(G)$ of a graph $G$ of maximum degree $\Delta(G)$
is either $\Delta(G)$ or $\Delta(G)+1$ \cite{vi},
Erd\H{o}s and Ne\v{s}et\v{r}il \cite{fashgytu2} posed the widely open conjecture $\chi_s'(G)\leq \frac{5}{4}\Delta(G)^2$.
Anderson \cite{an} and Hor\'{a}k, Qing, and Trotter \cite{hoqitr} 
showed $\chi_s'(G)\leq 10$ for graphs $G$ of maximum degree at most $3$.

Note that for a cubic graph $G$ of size $m(G)$, this last result implies 
$$\nu_s(G)\geq \frac{m(G)}{\chi_s'(G)}\geq \frac{m(G)}{10}.$$
In the present paper we strengthen this inequality further by proving the following.

\begin{theorem}\label{theorem1}
If $G$ is a cubic graph of size $m(G)$, then $\nu_s(G)\geq \frac{m(G)}{9}$.
\end{theorem}
Theorem \ref{theorem1} generalizes a recent result of Kang, Mnich, and M\"{u}ller, Corollary 3 in \cite{kamnmu},
who proved the same bound for cubic planar graphs.
Furthermore, it solves a conjecture posed by Henning and Rautenbach \cite{hera}.

Let the graph $K_{3,3}^+$ arise from $K_{3,3}$ by subdividing one edge once.
The connected graph $G$ that arises from the disjoint union of four copies of $K_{3,3}^+$ and two adjacent vertices $u$ and $v$
by joining each of $u$ and $v$ to two vertices of degree $2$ in the copies of $K_{3,3}^+$ 
is cubic and has $45$ edges. Since $\nu_s(G)=5$, this graph shows that Theorem \ref{theorem1} is best possible.
For a graph $G$, let $n_{3,3}^+(G)$ denote the number of components of $G$ that are isomorphic to $K_{3,3}^+$.
Furthermore, let $i(G)$ denote the number of isolated vertices of $G$.

Our Theorem \ref{theorem1} is a consequence of the following more detailed result.

\begin{theorem}\label{theorem2}
If $G$ is a subcubic graph of order $n(G)$, then 
$$\nu_s(G)\geq \frac{1}{6}(n(G)-i(G)-n_{3,3}^+(G)).$$
\end{theorem}
Theorem \ref{theorem2} actually implies the main result of Kang, Mnich, and M\"{u}ller, Theorem 1 in \cite{kamnmu},
stating that every subcubic planar graph $G$ of size $m(G)$ has an induced matching of size at least $m(G)/9$,
which can be determined in linear time.
In fact, since $K_{3,3}^+$ is not planar, every subcubic planar graph $G$ satisfies $n_{3,3}^+(G)=0$.
Furthermore, since $G$ is subcubic, we have $m(G)\leq \frac{3}{2}(n(G)-i(G))$ and hence
$$\nu_s(G)\geq \frac{1}{6}(n(G)-i(G)-n_{3,3}^+(G))\geq \frac{1}{9}m(G).$$
Our proof of Theorem \ref{theorem2} is constructive and 
leads to a linear time algorithm to determine an induced matching of the guaranteed size as we will explain below.
While the proof in \cite{kamnmu} is extremely involved, 
our proof of Theorem \ref{theorem2} is quite simple.
One reason for this simplicity is probably 
the fact that we express the lower bound on the strong matching number in terms of the order.
Since the strong matching number is the cardinality of an edge set, 
it might seem more natural to express this lower bound in terms of the size as in \cite{kamnmu}. 
Nevertheless, edges that are incident with low degree vertices should contribute more to the strong matching number
and our lower bound that is essentially linear in the order implicitly captures this intuitive idea.

In Section \ref{section2} we prove Theorem \ref{theorem2} and in 
Section \ref{section3} we conclude with some related results and problems.

\section{Proof of Theorem \ref{theorem2}}\label{section2}

For a contradiction, 
we assume that $G$ is a counterexample of minimum order,
that is, $\nu_s(G)<\frac{1}{6}(n(G)-i(G)-n_{3,3}^+(G))$.
The minimality of $G$ implies that $G$ is a connected graph of order at least $7$,
that is, $i(G)=n_{3,3}^+(G)=0$.

The {neighborhood} of a vertex $u$ in $G$ is denoted by $N_G(u)$
and the {closed neighborhood} $\{ u\}\cup N_G(u)$ of $u$ in $G$ is denoted by $N_G[u]$.
If $U$ is a set of vertices of $G$, 
then $G-U$ denotes the subgraph of $G$ induced by $V(G)\setminus U$.

\begin{claim}\label{c1}
$G$ does not contain $K_{3,3}^+$ as a subgraph.
\end{claim}
{\it Proof of Claim \ref{c1}:}
For a contradiction, 
we assume that $G$ has a subgraph $H$ that is isomorphic to $K_{3,3}^+$.
Let $u$ be the vertex of degree 2 in $H$.
Let $e$ be an edge of $H-N_H[u]$.
Let $G'$ be the graph obtained from $G$ by deleting the six vertices of degree 3 in $H$.
Since $G'$ is not an isolated vertex, we have $i(G')=0$.
Since every induced matching of $G'$ together with the edge $e$ is an induced matching of $G$,
we conclude
$\nu_s(G)
\geq \nu_s(G') +1
\geq \frac{1}{6}(n(G')-i(G'))+1
\geq \frac{1}{6}(n(G)-6)+1
= \frac{1}{6}n(G),$
which is a contradiction.
$\Box$

\medskip

\noindent A vertex of degree $1$ is an {\it end-vertex}.

\begin{claim}\label{c2}
The neighbor of an end-vertex has degree 3.
\end{claim}
{\it Proof of Claim \ref{c2}:}
For a contradiction, we assume that $u$ is an end-vertex of $G$ whose unique neighbor $v$ has degree $2$ in $G$.
Let $G'=G-N_G[v]$.
Clearly, $i(G')\leq 2$.
Since every induced matching of $G'$ together with the edge $uv$ is an induced matching of $G$,
we conclude
$\nu_s(G)
\geq \nu_s(G') +1
\geq \frac{1}{6}(n(G')-i(G'))+1
\geq \frac{1}{6}(n(G)-5)+1
>\frac{1}{6}n(G),$
which is a contradiction.
$\Box$

\begin{claim}\label{c3}
No two end-vertices have a common neighbor.
\end{claim}
{\it Proof of Claim \ref{c3}:}
For a contradiction, we assume that 
the two end-vertices $u_1$ and $u_2$ have the common neighbor $v$.
Let $G'=G-N_G[v]$.
Clearly, $i(G')\leq 2$.
Since every induced matching of $G'$ together with the edge $u_1v$ is an induced matching of $G$,
we conclude
$\nu_s(G)
\geq \nu_s(G') +1
\geq \frac{1}{6}(n(G')-i(G'))+1
\geq \frac{1}{6}(n(G)-6)+1
\geq \frac{1}{6}n(G),$
which is a contradiction.
$\Box$

\begin{claim}\label{c4}
No two end-vertices have distance $4$ in $G$.
\end{claim}
{\it Proof of Claim \ref{c4}:}
For a contradiction, we assume that 
$u_1v_1wv_2u_2$ is a shortest path in $G$
between the two end-vertices $u_1$ and $u_2$.
Let $G'=G-(N_G[v_1]\cup N_G[v_2])$.
By Claim \ref{c3}, we have $i(G')\leq 4$.
Since every induced matching of $G'$ together with the edges $u_1v_1$ and $u_2v_2$ is an induced matching of $G$,
we conclude
$\nu_s(G)
\geq \nu_s(G') +2
\geq \frac{1}{6}(n(G')-i(G'))+2
\geq \frac{1}{6}(n(G)-11)+2
>\frac{1}{6}n(G),$
which is a contradiction.
$\Box$

\begin{claim}\label{c5}
$\delta(G)\geq 2$.
\end{claim}
{\it Proof of Claim \ref{c5}:}
For a contradiction, we assume that $u$ is an end-vertex in $G$.
Let $v$ be its neighbor.
By Claim \ref{c2}, 
$v$ has degree 3.
We denote the two neighbors of $v$ distinct from $u$ by $w_1$ and $w_2$.
Let $G'=G-N_G[v]$.
By Claims \ref{c3} and \ref{c4},
we obtain $i(G')\leq 2$.
Since every induced matching of $G'$ together with the edge $uv$ is an induced matching of $G$,
we conclude
$\nu_s(G)
\geq \nu_s(G') +1
\geq \frac{1}{6}(n(G')-i(G'))+1
\geq \frac{1}{6}(n(G)-6)+1
=\frac{1}{6}n(G),$
which is a contradiction.
$\Box$

\begin{claim}\label{c6}
No two vertices of degree $2$ are adjacent.
\end{claim}
{\it Proof of Claim \ref{c6}:}
For a contradiction, we assume that $u_1$ and $u_2$ are two adjacent vertices of degree $2$.
Let $G'=G-(N_G[u_1]\cup N_G[u_2])$.
Note that $u_1$ and $u_2$ can have a common neighbor.
By Claim \ref{c5}, 
we obtain $i(G')\leq 2$.
Since every induced matching of $G'$ together with the edge $u_1u_2$ is an induced matching of $G$,
we conclude
$\nu_s(G)
\geq \nu_s(G') +1
\geq \frac{1}{6}(n(G')-i(G'))+1
\geq \frac{1}{6}(n(G)-6)+1
=\frac{1}{6}n(G),$
which is a contradiction.
$\Box$

\begin{claim}\label{c7}
No vertex of degree $2$ is contained in a triangle.
\end{claim}
{\it Proof of Claim \ref{c7}:}
For a contradiction, 
we assume that $u$ is a vertex of degree $2$ that has two adjacent neighbors $v_1$ and $v_2$.
By Claim \ref{c6}, 
the degrees of $v_1$ and $v_2$ are $3$.
Let $G'=G-N_G[v_1]$.
By Claim \ref{c5}, 
we obtain $i(G')\leq 1$.
Since every induced matching of $G'$ together with the edge $uv_1$ is an induced matching of $G$,
we conclude
$\nu_s(G)
\geq \nu_s(G') +1
\geq \frac{1}{6}(n(G')-i(G'))+1
\geq \frac{1}{6}(n(G)-5)+1
>\frac{1}{6}n(G),$
which is a contradiction.
$\Box$

\begin{claim}\label{c8}
No vertex of degree $2$ is contained in a cycle of length $4$.
\end{claim}
{\it Proof of Claim \ref{c8}:}
For a contradiction, 
we assume that $uv_1wv_2u$ is a cycle of length $4$ in $G$ 
such that $u$ is a vertex of degree $2$.
By Claims \ref{c6} and \ref{c7},
the vertices $v_1$ and $v_2$ are not adjacent and both of degree $3$.
Let $z$ be the neighbor of $v_1$ distinct from $u$ and $w$.
Let $G'=G-(N_G[v_1]\cup N_G[u])$.
Every induced matching of $G'$ together with the edge $uv_1$ is an induced matching of $G$.
Therefore, if $i(G')\leq 1$, then we obtain a similar contradiction as above.
This implies $i(G')\geq 2$.
By Claim \ref{c5} and 
since there are at most $4$ edges between $N_G[v_1]\cup N_G[u]$ and $V(G)\setminus (N_G[v_1]\cup N_G[u])$,
we obtain that $i(G')=2$ and that the two isolated vertices of $G'$, say $s_1$ and $s_2$, have degree $2$ in $G$.
This implies that $G$ is a uniquely determined graph of order $7$.
Furthermore, we may assume that $N_G(s_1)=\{ z,v_2\}$ and $N_G(s_2)=\{ z,w\}$.
Since $\{uv_2,s_2z\}$ is an induced matching of $G$,
we obtain $\nu_s(G)\geq 2>\frac{n(G)}{6}$,
which is a contradiction.
$\Box$

\begin{claim}\label{c9}
$G$ is cubic.
\end{claim}
{\it Proof of Claim \ref{c9}:}
For a contradiction, 
we assume that the vertex $u$ has degree $2$ in $G$.
Let $v_1$ and $v_2$ be the neighbors of $u$.
By Claims \ref{c6} and \ref{c7},
the vertices $v_1$ and $v_2$ are not adjacent and both of degree $3$.
Let $w_1$ and $w_2$ be the neighbors of $v_1$ distinct from $u$.
Let $s_1$ and $s_2$ be the neighbors of $v_2$ distinct from $u$.
By Claims \ref{c7} and \ref{c8} all seven vertices $u$, $v_1$, $v_2$, $w_1$, $w_2$, $s_1$, and $s_2$ are distinct.
Let $G'=G-(N_G[v_1]\cup N_G[u])$.
Every induced matching of $G'$ together with the edge $uv_1$ is an induced matching of $G$.
Therefore, if $i(G')\leq 1$, then we obtain a similar contradiction as above.
This implies $i(G')\geq 2$.
By Claims \ref{c5} and \ref{c8},
no isolated vertex $x$ of $G'$ satisfies $N_G(x)\subseteq \{ w_1,w_2\}$,
that is, every isolated vertex of $G'$ is adjacent to $v_2$.
This implies that $i(G')=2$ and that $s_1$ and $s_2$ are the two isolated vertices of $G'$.
Let $G''=G-(N_G[v_2]\cup N_G[u])$.
By symmetry, we obtain $i(G'')=2$ and that $w_1$ and $w_2$ are the two isolated vertices of $G''$.
This implies that $G$ has order $7$.
By Claim \ref{c1}, we may assume that $w_1$ and $s_2$ are not adjacent.
Now $\{ v_1w_1, v_2s_2\}$ is an induced matching of $G$,
we obtain $\nu_s(G)\geq 2>\frac{n(G)}{6}$,
which is a contradiction.
$\Box$

\medskip

\noindent Let $g(G)$ denote the length of a shortest cycle of $G$.

\begin{claim}\label{c10}
$g(G)>3$.
\end{claim}
{\it Proof of Claim \ref{c10}:}
For a contradiction, 
we assume that $v_1v_2v_3v_1$ is a triangle in $G$.
Let $G'=G-(N_G[v_1]\cup N_G[v_2])$.
By Claim \ref{c9} and since there are at most five edges between $N_G[v_1]\cup N_G[v_2]$
and $V(G)\setminus (N_G[v_1]\cup N_G[v_2])$,
we obtain $i(G')\leq 1$.
Since every induced matching of $G'$ together with the edge $v_1v_2$ is an induced matching of $G$,
we obtain a similar contradiction as above.
$\Box$

\begin{claim}\label{c11}
$g(G)>4$.
\end{claim}
{\it Proof of Claim \ref{c11}:}
For a contradiction, 
we assume that $v_1v_2v_3v_4v_1$ is a cycle of length $4$ in $G$.
By Claims \ref{c9} and \ref{c10},
the vertex $v_i$ has a neighbor $w_i$ outside of $\{ v_1,v_2,v_3,v_4\}$.
By Claim \ref{c10}, 
the vertices $w_1$ and $w_2$ are distinct.
Let $G'=G-(N_G[v_1]\cup N_G[v_2])$.
Every induced matching of $G'$ together with the edge $v_1v_2$ is an induced matching of $G$.
If $w_1$ and $w_2$ are adjacent,
then Claims \ref{c9} and \ref{c10} 
imply $i(G')=0$
and we obtain a similar contradiction as above.
Hence, by symmetry, 
we may assume that 
$w_1$ and $w_2$,
$w_1$ and $w_4$, as well as
$w_2$ and $w_3$ are not adjacent.
Now Claim \ref{c9} implies $i(G')=0$
and we obtain a similar contradiction as above.
$\Box$

\medskip

\noindent In view of the above claims, $G$ is cubic and $g(G)\geq 5$.
Let $u$ and $v$ be to adjacent vertices.
Let $G'=G-(N_G[u]\cup N_G[v])$.
Since $G$ is cubic and $g(G)\geq 5$,
we have $i(G')=0$.
Since every induced matching of $G'$ together with the edge $uv$ is an induced matching of $G$,
we obtain a similar contradiction as above,
which completes the proof of Theorem \ref{theorem2}.
$\Box$

\medskip

\noindent Before we proceed to the conclusion we comment on algorithmic aspects of proof of Theorem \ref{theorem2}.
Its eleven claims actually correspond to very simple local reductions.
Repeatedly applying these reductions to a given input graph in the order of their corresponding claims,
produces an induced matching in a recursive way.
Since $G$ is subcubic and each of these reductions only involves a bounded number of vertices,
the update of the remaining reduced graphs can be done in constant time.
Altogether, each edge of the induced matching found by this algorithm generates a constant effort,
which leads to an overall linear time algorithm.

\section{Conclusion}\label{section3}

As we have already pointed out Theorems \ref{theorem1} and \ref{theorem2} are best possible.
It seems an interesting problem to generalize them to larger maximum degrees, especially for regular graphs.

If $\Delta\geq 4$ is even and $G$ arises 
by replacing every vertex of a cycle of length $5$ with an independent set of order $\frac{\Delta}{2}$,
then $G$ is a $\Delta$-regular graph of order $\frac{5}{2}\Delta$ with $\nu_s(G)=1$.
Accordingly, we believe that for a regular graph $G$ of even maximum degree $\Delta(G)$, 
we have $\nu_s(G)\geq \frac{2n(G)}{5\Delta(G)}$,
which would also follow from the conjecture of Erd\H{o}s and Ne\v{s}et\v{r}il.
For odd regularity, one can generalize the construction explained after Theorem \ref{theorem1}.
Let $\Delta\geq 3$ be odd, say $\Delta=2r+1$.
Let $G_0$ arise by replacing the five vertices of a cycle of length $5$ 
with independent sets of order $r+1$, $r+1$, $r$, $r$, and $r$ in cyclic order, respectively. 
Note that $G_0$ has $2r+1$ vertices of degree $\Delta$ and $r$ vertices of degree $\Delta-1$.
Let $G$ arise from the disjoint union of four copies of $G_0$, say $G_0^1$, $G_0^2$, $G_0^3$, and $G_0^4$, 
and two adjacent vertices $u$ and $v$ by joining 
$u$ to the vertices of degree $\Delta-1$ in $G_0^1$ and $G_0^2$
and 
$v$ to the vertices of degree $\Delta-1$ in $G_0^3$ and $G_0^4$.
Clearly, $G$ is a $\Delta$-regular graph of order $10\Delta$ with $\nu_s(G)=5$.
Accordingly, we believe that for a regular graph $G$ of odd maximum degree $\Delta(G)$, 
we have $\nu_s(G)\geq \frac{n(G)}{2\Delta(G)}$,
which is better than the immediate consequence of the conjecture of Erd\H{o}s and Ne\v{s}et\v{r}il.

A very simple greedy argument shows that 
if $T$ is a forest of maximum degree $\Delta(T)$, 
then $\nu_s(T)\geq \frac{m(T)}{2\Delta(T)-1}$.
Similarly, if $G$ is a graph of maximum degree $\Delta(G)$, 
then $\nu_s(G)\geq \frac{m(G)}{2\Delta(G)(\Delta(G)-1)+1}$.
As for bounds in terms of the order, 
a first result along the lines of the proof of Theorem \ref{theorem2} is the following.

\begin{proposition}\label{proposition1}
If $G$ is a graph of girth at least $6$, then 
$\nu_s(G)\geq \frac{n(G)-i(G)}{\frac{1}{4}\Delta(G)^2+\Delta(G)+1}$.
\end{proposition}
{\it Proof:}
For a contradiction, we assume that $G$ is a counterexample of minimum order.
Clearly, $G$ is connected, $i(G)=0$, and $n(G)>\frac{1}{4}\Delta(G)^2+\Delta(G)+1$.

If $G$ has no end-vertex, then let $uv$ be an edge of $G$.
Let $G'=G-(N_G[u]\cup N_G[v])$.
Since $G$ has no cycle of length at most $5$, we have $i(G')=0$.
Since every induced matching of $G'$ together with the edge $uv$ is an induced matching of $G$,
we obtain 
$\nu_s(G)
\geq \nu_s(G')+1
\geq  \frac{n(G')}{\frac{1}{4}\Delta(G)^2+\Delta(G)+1}+1
\geq \frac{n(G)-2\Delta(G)}{\frac{1}{4}\Delta(G)^2+\Delta(G)+1}+1
\geq \frac{n(G)}{\frac{1}{4}\Delta(G)^2+\Delta(G)+1}$.
Hence we may assume that $G$ has end-vertices.

Let $k$ denote the maximum cardinality of a set of end-vertices that have a common neighbor.
Let the vertex $v$ be adjacent to $k$ end-vertices.
Let $u$ be an end-vertex that is a neighbor of $v$.
Let $G'=G-N_G[v]$.
By the definition of $k$,
we have $n(G)-n(G')=d_G(v)+1\leq \Delta(G)+1$
and $i(G')=k(d_G(v)-k)\leq k(\Delta(G)-k)$.
This implies
$n(G)-(n(G')-i(G'))\leq \Delta(G)+1+k(\Delta(G)-k)\leq \frac{1}{4}\Delta(G)^2+\Delta(G)+1$.
Since every induced matching of $G'$ together with the edge $uv$ is an induced matching of $G$,
we obtain 
$\nu_s(G) \geq \nu_s(G')+1 \geq \frac{n(G)}{\frac{1}{4}\Delta(G)^2+\Delta(G)+1},$
which completes the proof.
$\Box$

\medskip

\noindent A best possible version of Proposition \ref{proposition1} would be interesting.


\begin{thebibliography}{0}
\bibitem{an}
L.D. Andersen,
The strong chromatic index of a cubic graph is at most 10,
Discrete Math. 108 (1992) 231-252.

\bibitem{brmo}
A. Brandst\"{a}dt and R. Mosca,
On distance-3 matchings and induced matchings,
Discrete Appl. Math. 159 (2011) 509-520.

\bibitem{ca1}
K. Cameron,
Induced matchings,
Discrete Appl. Math. 24 (1989) 97-102.

\bibitem{ca2}
K. Cameron,
Induced matchings in intersection graphs,
Discrete Math. 278 (2004) 1-9.

\bibitem{ed}
J. Edmonds,
Paths, trees, and flowers,
Canad. J. Math. 17 (1965) 449-467.

\bibitem{fashgytu2}
R.J. Faudree, R.H. Schelp, A. Gy\'{a}rf\'{a}s, and  Zs. Tuza,
The strong chromatic index of graphs,
Ars Comb. 29B (1990) 205-211.

\bibitem{gohehela}
W. Goddard, S.M. Hedetniemi, S.T. Hedetniemi, and R. Laskar,
Generalized subgraph-restricted matchings in graphs,
Discrete Math. 293 (2005) 129-138.

\bibitem{hera}
M.A. Henning and D. Rautenbach,
Induced matchings in subcubic graphs without short cycles,
to appear in Discrete Math.

\bibitem{hoqitr}
P. Hor\'{a}k, H. Qing, and W.T. Trotter,
Induced Matchings in Cubic Graphs,
J. Graph Theory 17 (1993) 151-160.

\bibitem{kamnmu}
R.J. Kang, M. Mnich, and T. M\"{u}ller,
Induced matchings in subcubic planar graphs,
SIAM J. Discrete Math. 26 (2012) 1383-1411.

\bibitem{lopl}
L. Lov\'{a}sz and M.D. Plummer,
Matching Theory,
vol. 29, Annals of Discrete Mathematics, North-Holland, Amsterdam, 1986.

\bibitem{lo}
V.V. Lozin,
On maximum induced matchings in bipartite graphs,
Inf. Process. Lett. 81 (2002) 7-11.

\bibitem{stva}
L.J. Stockmeyer and V.V. Vazirani,
NP-completeness of some generalizations of the maximum matching problem,
Inf. Process. Lett. 15 (1982) 14-19.

\bibitem{tu}
W.T. Tutte,
The factorization of linear graphs,
J. Lond. Math. Soc. 22 (1947) 107-111.

\bibitem{vi}
V.G. Vizing,
On an estimate of the chromatic class of a p-graph,
Diskret. Analiz. 3 (1964) 25-30.

\end{thebibliography}
\end{document}